\documentclass[a4paper,twoside,12pt]{article}
\usepackage[english]{babel}
\usepackage{bbding}
\usepackage[latin1]{inputenc}
\usepackage{kpfonts}
\usepackage{fancyhdr}
\usepackage{amsmath,amssymb,color,epsfig}
\usepackage{mathrsfs}
\usepackage{eufrak}
\usepackage{dsfont}
\usepackage{upgreek}
\usepackage{fancyhdr}

\newtheorem{theorem}{Theorem}[section]
\newtheorem{proposition}[theorem]{Proposition}

\newtheorem{lemma}[theorem]{Lemma}

\def\tr{{\rm tr}}

\def\vp{{\varphi}}

\def\1{\mathds{1}}
\def\e{\varepsilon}

\title{Generalized $\beta$-Gaussian Ensemble Equilibrium measure method}
\date{}
\begin{document}
\author{Mohamed BOUALI}
\maketitle
\begin{center}{\bf Abstact}
   \end{center}
We describe $\beta$-Generalized random Hermitian matrices ensemble sometimes called Chiral ensemble.  We give global asymptotic of the density of eigenvalues or the statistical density. We investigate general method names as equilibrium measure method. When taking $n$ large limit we will see that the asymptotic density of eigenvalues generalize the Wigner semi-circle law.
\section{Introduction} The generalized $\beta$-Gaussian ensemble, generalize the classical random matrix ensemble: Gaussian orthogonal,
unitary and symplectic ensembles (denoted by GOE, GUE and GSE for short, which correspond to
the Dyson index $\beta$ = $1$, $2$ and $4$), from the quantization index to the continuous exponents $\beta > 0$. These
ensembles possess the joint probability density function (p.d.f.) of real eigenvalues $\lambda_1, . . . , \lambda_n$ with the
form
$$\mathbb{P}_n(d\lambda)=\frac{1}{Z_n}e^{-\sum\limits_{i=1}^n\lambda_i^2}\prod_{1\leq i<j\leq n}|\lambda_i-\lambda_j|^{\beta}d\lambda_1...d\lambda_n,$$
where $Z_n$ can be evaluated by the using the Selberg integral
$$Z_n=(2\pi)^{2n}\prod_{i=1}^n\frac{\Gamma(1+\frac{\beta i}{2})}{\Gamma(1+\frac{\beta}{2})}.$$
Recently, Dumitriu and Eldeman have construct a tri-diagonal matrix model of these ensembles see \cite{I}.

Basing on the p.d.f. of eigenvalues $\mathbb{P}_n$, the (level) density, or one-dimensional marginal eigenvalue
density scaled by the factor $\frac1{\sqrt {2n}}$ converge weakly to the famous Wigner semi-cercle law as follows:
for every bounded continuous functions $f$ on $\Bbb R$
$$\lim_{n\to\infty}\int_{\Bbb R}f(\frac{x}{\sqrt{2n}})h_n(x)dx=\int_{\Bbb R}f(x)\rho(x)dx,$$
where $$\rho(x)=\left\{\begin{aligned}&\frac{1}{\pi}\sqrt{2-x^2}\;\;{\rm if}\;|x|\leq\sqrt2\\&0\quad\qquad\quad\;\;{\rm if}\;|x|\geq\sqrt2 \end{aligned}\right.$$
$$h_n(\lambda_1)=\int_{\Bbb R^{n-1}}\mathbb{P}_n(\lambda_1,...,\lambda_n)d\lambda_2...d\lambda_n.$$
Many others work in this direction of random matrices and asymptotic of eigenvalues has been developed in the last years, one can see \cite{fo}, for a good reference.

In this work we will study a generalization of the Gaussian random matrices ensemble which is called some times the Chiral-ensemble when $\beta=1, 2$ or $4$. We will consider the general case where $\beta>0$, in that case the joint probability density in $\Bbb R^n$ is given by:
$$\mathbb{P}_n(dx)=\frac{1}{Z_n}e^{-\sum\limits_{i=1}^nx_i^2}\prod_{i=1}^n|x_i|^{2\lambda}\prod_{1\leq i<j\leq n}|x_i-x_j|^{\beta}dx_1...dx_n,$$
where $Z_n$ is a normalizing constant and $\lambda$ is a positive parameter.
Using a general method of logarithmic potential we will prove that, the statistical density of eigenvalues converge for the tight topology as $n\to +\infty$ to some probability density. Which generalized the Wigner semi-circle law. Such result has been proved in \cite{I0} for $\beta=2$, by the orthogonal polynomials method.

The paper is organized as follow.
In sections 2 and 3 we gives some results about classical potential theory, which will be used together with some fact about boundary values distribution to characterized the Cauchy transform of some equilibrium measures.

In section 4, we will describe the model to study, as physics model, and we give the joint probability density. Moreover we defined the statistical density $\nu_n$ of eigenvalues and we explain how the eigenvalues must be rescaling by the factor $\sqrt n$. Also we gives the first means result theorem 4.1, which state the convergence of the statistical density $\nu_n$ to some probability measure $\nu_{\beta, c}$. We will prove that, the measure $\nu_{\beta, c}$ is an equilibrium measure and we compute the exact value of the energy for general $\beta$, after calculating the energy for $\beta=2$.

In section 5  we gives the proof of the first result of theorem 4.1.

\section{Logarithmic potential}
The logarithmic potential of a positive measure $\nu$ on $\Bbb R$ is the function $U^{\nu}$ defined by
$$U^{\nu}(x)=\int_{\Bbb R}\log{\frac{1}{|x-t|}}\nu(dt).$$
It will defined with value on $]-\infty,+\infty]$ if $\nu$ is with compactly support or more general, if
$$\int_{\Bbb R}\log{(1+|t|)}\nu(dt)<\infty.$$
Observe that $$\lim_{n\to\infty}(U^{\nu}(x)+\nu(\Bbb R)\log|x|)=0.$$
The Cauchy transform $G_{\nu}$ of a bounded measure $\nu$ on $\Bbb R$ is the function
defined on  ${\Bbb C}\setminus {\rm supp}(\nu)$ by
$$G_{\nu}(z)=\int_{\Bbb R}\frac{1}{z-t}\nu(dt).$$
The Cauchy transform is holomorphic.\\
Assume that ${\rm supp}(\nu) \subset] -\infty, a]$, and
$$\int_{\Bbb R}\log{(1+|t|)}\nu(dt)<\infty.$$
Then the function
$$F(z)=\int_{\Bbb R}\log{\frac{1}{|z-t|}}\nu(dt).$$
is defined and holomorphic in $\Bbb C\setminus]-\infty,a]$. Furthermore $\displaystyle F'(z)=G_{\nu}(z),$ and
$$\begin{aligned}
&U^{\nu}(x)=-{\rm Re}F(x)\qquad(x>a)\\
&U^{\nu}(x)=-\lim_{\e\to 0}{\rm Re}F(x+i\e)\qquad(x\in\Bbb R)
\end{aligned}$$
In the distribution sense,
$$\frac{d}{dx}U^{\nu}(x)=-{\rm Re}G_{\nu}(x).$$
We will use some properties of the boundary value distribution of a
holomorphic function. Let $f$ be holomorphic in $\Bbb C\setminus\Bbb R$. It is said to be of
moderate growth near $\Bbb R$ if, for every compact set $K \subset\Bbb R$, there are $\e > 0$, $N>0$, and $C>0$ such that
$$|f(x+iy)|\leq\frac{C}{|y|^N}\quad (x\in K, 0<|y|\leq\e).$$
Then for all $\varphi\in{\cal D}{(\Bbb R)}$,
$$(T,\varphi)=\lim_{\e\to 0, \e>0}\int_{\Bbb R}\vp(x)(f(x+i\e)-f(x-i\e))dx,$$
defines a distribution $T$ on $\Bbb R$. It is denoted $T = [f]$, and called the
difference of boundary values of $f$. One shows that the function $f$ extends
as a holomorphic function in $\Bbb C\setminus {\rm supp}([f])$. In particular, if $[f] = 0$, then
$f$ extends as a holomorphic function in $\Bbb C$.\\
For $\alpha\in\Bbb C$, the distribution $Y_{\alpha}$ is defined, for ${\rm Re}\alpha > 0$, by
$$(Y_{\alpha},\vp)=\frac{1}{\Gamma(\alpha)}\int_{\Bbb R}\vp(t)t^{\alpha-1}dt.$$
The distribution $Y_{\alpha}$, as a function of $\alpha$, admits an analytic continuation
for $\alpha\in\Bbb C$. In particular $Y_0 = \delta_0$, the Dirac measure at $0$.

For $\alpha\in\Bbb C$, we defines the holomorphic function $z^{\alpha}$ in $\Bbb C\setminus\,] - \infty, 0]$ as
follows: if $z = re^{i\theta}$, with $r > 0$, $-\pi < \theta <\pi$, then

$$z^{\alpha}=r^{\alpha}e^{i\alpha\theta}.$$
The function $z^\alpha$ is of moderate growth near $\Bbb R$, and

$$([z^{\alpha}],\vp)=-2i\pi\frac{1}{\Gamma(-\alpha)}(Y_{\alpha+1},\check\vp),$$
where $\check\vp(t)=\vp(-t)$. In particular when $\alpha=-1$
$$[\frac 1z]=-2i\pi\delta_0.$$

\begin{proposition}
Let $\nu$ be a bounded positive measure on $\Bbb R$.\\
$(\rm i)$ The Cauchy transform $G_{\nu}$ of $\nu$ is holomorphic in $\Bbb C\setminus{\rm supp}(\nu)$, of
moderate growth near $\Bbb R$, and
$$[G_{\nu}]=-2i\pi\nu.$$
$(\rm ii)$ Assume that the support of $\nu$ is compact. Let $F$ be holomorphic in
$\Bbb C \setminus\Bbb R$, of moderate growth near $\Bbb R$, such that
$$[F]=-2i\pi\nu.$$
Then F is holomorphic in $\Bbb C\setminus{\rm supp(\nu)}$. If further
$$\lim_{|z|\to\infty}F(z)=0.$$
Then $$G_{\nu}=F.$$

\end{proposition}
\section{ Equilibrium measure some basic results}
Let us first recall some basic
facts about the tight topology. All the present result in equilibrium measure can be find in the good reference \cite{sa} and references therein. Let ${\mathfrak M}^1(\Sigma)$ be the set of probability
measures on the closed set $\Sigma\subset\Bbb R$. We consider the tight topology. For this
topology a sequence $(\nu_n)$ converges to a measure $\nu$ if, for every continuous
bounded function $f$ on $\Sigma$,
$$\lim_{n\to\infty}\int_{\Sigma}f(x)\nu_n(dx)=\int_{\Sigma}f(x)\nu(dx).$$
This topology is metrizable. If $\Sigma$ is bounded, then ${\mathfrak M}^1(\Sigma)$ is compact.\\
Let $\Sigma$ be a closed interval $(\Sigma = \Bbb R, ]-\infty,a], [ b, +\infty [ or [a, b])$, and $Q$ a
function defined on $\Sigma$ with values on $] -\infty,+\infty]$, continuous on int$(\Sigma)$. If
$\Sigma$ is unbounded, it is assumed that
$$\lim_{|x|\to +\infty}(Q(x)-\log(1+x^2))=+\infty.$$
If $\nu$ is a probability measure supported by $\Sigma$, the energy $E(\nu)$ of $\nu$ is
defined by
$$E(\nu)=\int_{\Sigma\times\Sigma}\log\frac{1}{|x-y|}\nu(dx)\nu(dy)+\int_{\Sigma}Q(x)\nu(dx).$$
which mean that
$$E(\nu)=\int_{\Sigma}U^{\nu}(x)\nu(dx)+\int_{\Sigma}Q(x)\nu(dx).$$
 By a straightforward computation we can prove that $E(\nu)$ is bounded below. Hence we defined
 $$\displaystyle E^*=\inf\{E(\nu)\mid \nu\in\mathfrak{M}^1(\Sigma)\}.$$
\begin{theorem}
If $\nu(dx)=f(x)dx$, where $f$ is a continuous function with
compact support $\subset int(\Sigma)$. Then the potential $U^{\nu}$ is a continuous function, and
$E^*\leq E(\nu) < \infty$. Furthermore there is a unique measure $\nu^*\in\mathfrak{M}^1(\Sigma)$ such that
$$E^*=E(\nu^*).$$
The support of $\nu^*$ is compact.\\
This measure $\nu^*$ is called the equilibrium measure.
\end{theorem}
\begin{proposition}
Let $\nu\in\mathfrak{ M}^1(\Sigma)$ with compact support.
Assume that the potentiel $U^{\nu}$ of $\nu$ is continuous and that there is a
constant $C$ such that\\
${\rm (i)}$ $U^{\nu}(x)+\frac 12Q(x)\geq C$ on $\Sigma$.\\
${\rm (ii)}$ $U^{\nu}(x)+\frac 12Q(x)= C$ on ${\rm supp}(\nu)$. Then $\nu$ is the equilibrium measure: $\nu=\nu^*$.
\end{proposition}
The constant $C$ is called the (modified) Robin constant. Observe that
$$E^*=C+\frac 12\int_{\Sigma}Q(x)\nu^*(dx).$$
It is easy to see the action by linear transformation on the energy.
\begin{proposition}
Let the transformation $h(s) = as + b$ map $\Sigma$ onto
$\Sigma'$. If $Q$ is defined on $\Sigma'$, then $Q\circ h$ is defined on $\Sigma$. If $\nu$ is a probability
measure on $\Sigma$, then $\sigma = h(\nu)$ is the probability measure on $\Sigma'$ defined
by $$\int_{\Sigma'}f(t)\sigma(dt)=\int_{\Sigma}f\circ h(t)\nu(dt).$$
Then
$$E_{(\Sigma', Q)}(h(\nu))=E_{(\Sigma, Q\circ h)}(\nu)-\log|a|.$$
\end{proposition}
For the proof of the previous theorem and proposition, see for instance theorem II.2.3, proposition II.3.1 of \cite{F}.

\section{Statistical of the generalized Gaussian unitary ensemble}

Let $H_n=Herm(n, \Bbb F)$ be the vector space of square Hermitian matrices with coefficient in the field $\Bbb F=\Bbb R,\,\Bbb C$ or $\Bbb H$. For $\mu>-\frac 12$, we denote by $\mathbb{P}_{n,\mu}$
the probability measure on $H_n$ defined by.
$$\int_{H_n}f(x)\mathbb{P}_{n, \mu}(dx)=\frac{1}{C_n}\int_{H_n}f(x)|\det(x)|^{2\mu}e^{-\tr(x^2)}m_n(dx),$$
for a bounded mesurable function $f$, where $m_n$ is the Euclidean measure associated to the usual inner product $<x, y>=\tr(xy)$ on $H_n$ and $C_n$ is a normalized constant. which is given for $d=2$ by
$$C_n=n!\prod_{k=0}^{n-1}\gamma_\mu(k)$$
\begin{equation}\gamma_\mu(k)=\left\{\begin{aligned}&m!\Gamma(m+\mu+\frac12)\quad\mbox {if}\; k=2m,\\
&m!\Gamma(m+\mu+\frac32)\quad\mbox {if}\; k=2m+1.\\
\end{aligned}\right.\end{equation}
For general $\beta=1$ or $4$ the constant is given by Jack polynomials.

When $\mu=0$ we recover the classical Gaussian unitary ensemble and, $$\displaystyle C_n=\pi^{\frac n2}2^{\frac{n(n-1)}{2}}\prod_{k=0}^{n}k!.$$
We endowed the space $H_n$ with the probability measure $\mathbb{P}_{n, \mu}$. The probability $\mathbb{P}_{n,\mu}$ is invariant for the action of the unitary group $U(n)$ by the conjugation
$$x\mapsto uxu^*\qquad (u\in U(n)).$$
\subsection{Spectral density of eigenvalues}

Let $f$ be a $U(n)$-invariant function on $H_n$.
$$f(uxu^*)=f(x)\qquad\forall\;u\in U(n),$$
Then by the spectral theorem there exist a symmetric function $F$ in $\Bbb R^n$ such that
$$f(x)=F(\lambda_1,...,\lambda_n).$$
If $f$ is integrable with respect to $\mathbb{P}_{n,\mu}$, then by using the formula of integration of Well we obtain
$$\int_{H_n}f(x)\mathbb{P}_{n,\mu}(dx)=\int_{\Bbb \Bbb^n}F(\lambda_1,\ldots,\lambda_n)q_{n, \mu}(\lambda_1,...,\lambda_n)d\lambda_1\ldots\lambda_n,$$
where
$$q_{n,\mu}(\lambda_1,\ldots,\lambda_n)=\frac{1}{c_n}e^{-\sum\limits_{k=1}^{n}\lambda_k^2}\prod\limits_{k=1}^{n}|\lambda_k|^{2\mu}
|\Delta(\lambda)|^{\beta},$$
and $\beta=1,\,2,\,4$ for $\Bbb F=\Bbb R$, $\Bbb C$ or $\Bbb H$. $\Delta$ is the vandermonde determinant

$$\Delta(\lambda)=\prod_{1\leq i<j\leq n}(\lambda_i-\lambda_j),$$ .\\

More general we will consider $n$ particles free to move in $\Bbb R^n$, in equilibrium at absolute temperature $T$. A fundamental postulate gives the
p.d.f. for the event that the particles are at positions $\lambda_1, . . . ,\lambda_n$ as:
$$q_{n,\mu_n}(\lambda_1,\ldots,\lambda_n)=\frac{1}{Z_N}e^{-\beta V_n(\lambda_1,...\lambda_n)},$$

where $$V_n(\lambda_1,...,\lambda_n)=\frac{2n}{\beta}\sum\limits_{k=1}^{n}(\lambda_k^2+\frac{2\mu_n}{n}\log\frac{1}{|\lambda_k|})+
\sum_{1\leq i<j\leq n}\log\frac{1}{|\lambda_i-\lambda_j|}.$$
Here $V_n(\lambda_1,...,\lambda_n)$ denotes the total potential energy of the system, $\beta:= \frac{1}{k_BT }$ ($k_B$ is Boltzmann\'s
constant), and $Z_n$ is a normalizing constant.

The term $V_n(\lambda_1,...,\lambda_n)$ is referred to as the Boltzmann factor and  ${\widetilde Z}_n:=\frac{Z_n}{n!} $ is called the (canonical) partition function.

Our first result is to study as n go to infinity the asymptotic of the Normalized Counting
Measure (Density of States) $\nu_n$ defined on $\Bbb R$ as follows: if $f$ is a measurable function,
$$\int_{\Bbb R}f(t)\nu_n(dt)=\mathbb{E}_{n,\mu_n}(\frac{1}{n}\sum_{i=1}^nf(\lambda_i)),$$
where  $\mathbb{E}_{n,\mu_n}$ is the expectation with respect the probability measure on $\Bbb R^n$
$$\mathbb{P}_{n,\mu_n}(d\lambda)=\frac{1}{Z_n}e^{-n\sum\limits_{k=1}^{n}Q_n(\lambda_k)}
\prod_{1\leq i<j\leq n}|\lambda_i-\lambda_j|^{\beta},$$
and $$Q_n(x)=x^2+\frac{2\mu_n}{n}\log\frac{1}{|x|}.$$
By invariance of the measure $\mathbb{P}_{n,\mu_n}$ by the symmetric group, we have that the
 measure $\nu_n$ is continue with respect to the Lebesgue measure
$$\nu_n(dt)=h_{n,\mu_n}(t)dt.$$
where $$h_{n,\mu_n}(t)=\int_{\Bbb R^{n-1}}q_{n,\mu_n}(t,\lambda_2,\ldots,\lambda_n)d\lambda_2\ldots d\lambda_n.$$\\

Let compute the two first moments of the measure $\nu_n$:
$$m_1(\mu_n)=\int_{\Bbb R}t\nu_n(dt)=\frac 1n\int_{\Bbb R^n}\sum_{k=1}^{n}\lambda_k\;\mathbb{P}_{n,\mu_n}(d\lambda)=0,$$
the second moment is:
$$m_2(\nu_n)=\frac 1n\int_{\Bbb R^n}\sum_{k=1}^{n}\lambda_k^2\;\mathbb{P}_{n,\mu_n}(d\lambda),$$
Since for all $\alpha>0$, $$Z_n(\alpha)=\int_{\Bbb R^n}e^{-\alpha\sum\limits_{k=1}^{n}\lambda_k^2}\prod\limits_{k=1}^{n}|\lambda_k|^{2\mu_n}
\prod_{1\leq i<j\leq n}|\lambda_i-\lambda_j|^{\beta}d\lambda_1\cdots d\lambda_n=\alpha^{-n\mu_n-\frac{\beta}{4}n(n-1)-\frac n2}Z_n,$$
and $$m_2(\nu_n)=-\frac 1n\frac{d}{d\alpha}\log(Z_n(\alpha))|_{\alpha=1}=\mu_n+\frac{\beta}{4}(n-1)+\frac12.$$

This suggests that $\nu_n$ does not converge, and that a scaling of order $\displaystyle\sqrt {n\frac{\beta}{4}+\mu_n}$ is necessary.

 We come to The mean result: the measure $\nu_n$ converge weakly to some probability measure $\nu_{\beta, c}$ which is an equilibrium measure.

\begin{theorem}
Let $(\mu_n)_n$ be a nonnegative real sequence, if $\displaystyle\lim_{n\to\infty}\frac{\mu_n}{n}=c$. Then
the probability measure $\nu_n$ converge weakly to the probability $\nu_{\beta, c}$, where $\nu_{\beta, c}$ is the measure on $S=[-b,-a]\cup[a,b]$ with density with respect to the Lebesque measure $$f_{\beta,c}(t)=\left\{\begin{aligned}&\frac {2}{\pi\beta}\frac{1}{|t|}\sqrt{(t^2-a^2)(b^2-t^2)}\;\;\;if\;\;t\in S\\
&\quad 0\qquad\qquad\qquad\qquad\qquad\;\; if\;\;t\notin S\end{aligned}\right.,$$
and $a=\sqrt{\frac{\beta}{2}}\sqrt{1+\frac{2c}{\beta}-\sqrt {1+\frac{4c}{\beta}}}$, $b=\sqrt{\frac{\beta}{2}}\sqrt{1+\frac{2c}{\beta}+\sqrt {1+\frac{4c}{\beta}}}$.
Moreover the energy of the equilibrium measure $\nu_{\beta,c}$ is
$$ E^*_{\beta,c}=\frac{3\beta}{8}+\frac{\beta}{4}\log(\frac{4}{\beta})+c(\frac{3}{2}+\log\frac4{\beta})+\frac{2c^2}\beta\log\frac{4c}{\beta}
-(\frac{2c^2}{\beta}+c+\frac\beta8)\log(1+\frac{4c}{\beta}).$$
\end{theorem}
The convergence is in the sense that for every continuous bounded function $f$ on $\Bbb R$
$$\lim_{n\to\infty}\int_{\Bbb R}f(t)\nu_n(dt)=\int_{\Bbb R}f(t)\nu_{\beta, c}(dt).$$

\subsection{Equilibrium measure of generalized Gaussian unitary ensemble}

For $c\geq 0$, $\beta>0$, one considers on $\Sigma=\Bbb R$, the potential
$$Q_{c}(t)=t^2+2c\log\frac {1}{|t|},$$

The energy of a probability measure $\mu\in\mathfrak{M}^1(\Bbb R)$ is defined by
$$E_{\beta, c}(\mu)=\frac{\beta}{2}\int_{\Bbb R^2}\log\frac{1}{|s-t|}\mu(ds)\mu(dt)+\int_{\Bbb R}Q_c(t)\mu(dt),$$
and let $U_{\beta, c}$ be the potential of the measure $\nu_{\beta, c}$, $\displaystyle U_{\beta, c}(x)=\int_{\Bbb R}\log\frac{1}{|x-y|}\nu_{\beta, c}(dy)$
\begin{proposition} The probability measure $\nu_{\beta, c}$ is the equilibrium measure, which mean that
$$\inf\left\{E_{\beta, c}(\nu)\mid \nu\in\mathfrak{M}^1(\Bbb R)\right\}=E(\nu_{\beta, c})=E^*_{\beta,c}.$$
Furthermore \\
$(i)$ $\displaystyle U_{\beta, c}(x)+\frac12 Q_c(x)= C, $ on $S$.\\
$(ii)$ $\displaystyle U_{\beta, c}(x)+\frac12 Q_c(x)\geq C,$ on $\Bbb R\setminus S$.\\

\end{proposition}
We will give the value of the energy $E^*_{\beta, c}$ in section 3.3.\\
To prove the proposition we need same preliminary results and then applying proposition 2.2. \\
For more convenient notation we shall denote $\displaystyle c'=\frac{2c}{\beta}$.\\

Putting $$f(z)=\frac 2{\beta z}\sqrt{z-a}\sqrt{z-b}\sqrt{z+a}\sqrt{z+b},$$
The function $f$ is holomorphic on the domain $\Bbb C\setminus (S\cup\{0\})$, of
moderate growth near $S\cup\{0\}$.

\begin{proposition}The difference between the two limits values of $f$ in the distribution sense,
 $[f]=f(x+i0)-f(x-i0),$ is given by
 $$[f]=2i\pi\nu_{\beta, c'}+2i\pi c'\delta_0.$$
 \end{proposition}
 {\bf Proof.}\\
For $b>0$, observe that the function $\displaystyle f(z)=\frac 2{\beta z}\sqrt{z-a}\sqrt{z-b}\sqrt{z+a}\sqrt{z+b},$ is defined and holomorphic on $\Bbb C\setminus ]-\infty,b].$\\
 For $x>b$, $\displaystyle f(x)=\frac 2{\beta x}\sqrt{(x-a)(x-b)(x+a)(x+b)}=\frac 2{\beta x}\sqrt{(x^2-a^2)(x^2-b^2)},$ be the usual square root of positive numbers.\\
 For $x<-b$, $$\lim_{\e\to 0,\e>0}f(x\pm i\e)=e^{\pm i\pi}e^{\pm i\pi}\frac 2{\beta x}\sqrt{(a-x)(b-x)(-x-a)(-x-b)}=\frac 2{\beta x}\sqrt{(x^2-a^2)(x^2-b^2)}.$$
 There for $f$ extended as holomorphic function on $\Bbb C\setminus S$. Furthermore
 For $-b<x<-a$, $$\lim_{\e\to 0,\e>0}f(x\pm i\e)=-e^{\pm i\frac{\pi}{2}}\frac 2{\beta x}\sqrt{(a-x)(b-x)(-x-a)(x+b)}=\mp i\frac 2{\beta x}\sqrt{(a^2-x^2)(b^2-x^2)}.$$
 For $-a<x<a$, $x\neq 0$, $$\lim_{\e\to 0,\e>0}f(x\pm i\e)=e^{\pm i\pi}\frac 2{\beta x}\sqrt{(a-x)(b-x)(x+a)(x+b)}=-\frac 2{\beta x}\sqrt{(a^2-x^2)(b^2-x^2)}.$$
 For $z$ near $0$, by Taylor expansion $$f(z)=-\frac{2ab}{\beta z}+g(z)=\frac{-c'}{z}+g(z),$$
 where $g$ is an holomorphic function.\\

 For $a<x<b$, $$\lim_{\e\to 0,\e>0}f(x\pm i\e)=e^{\pm i\frac{\pi}{2}}\frac 2{\beta x}\sqrt{(a-x)(b-x)(-x-a)(x+b)}=\pm i\frac 2{\beta x}\sqrt{(x^2-a^2)(b^2-x^2)}.$$
   It follows that
   $$[f]=2i\frac 2{\beta |x|}\sqrt{(x^2-a^2)(b^2-x^2)}\,\chi_S+2i\pi c'Y_0=2i\pi\nu_{\beta, c'}+2i\pi c'\delta_0.$$
   Which complete the proof of the proposition.\qquad\qquad\qquad\qquad\qquad\qquad\qquad\qquad\qquad$\blacksquare$ \\

Let denote by $G_{\beta, c'}$ the Cauchy transform of the measure $\nu_{\beta, c'}$: for all $z\in\Bbb C\setminus S $,
 $$G_{\beta, c'}(z)=\int_\Bbb R\frac{1}{z-t}\nu_{\beta, c'}(dt).$$

 \begin{proposition}
 The Cauchy transform of the measure $\nu_{\beta, c'}$ is defined on $\Bbb C\setminus S $,
 $$G_{\beta, c'}(z)=-f(z)+\frac{2}{\beta}z-\frac {c'}z$$
 \end{proposition}
 {\bf Proof.}
 From the previous proposition we have for all $x\in S$
 $$\lim_{\e\to 0, \e>0}(f(x+i\e)-f(x-i\e))=2i\frac 2{\beta x}\sqrt{(x^2-a^2)(b^2-x^2)}.$$
 It follows that, if $\vp$ is a holomorphic function in a neighborhood $U$ of
$S$, and $\gamma$ is a path in $U$ around $S$ in the positive sense, then
$$\frac{1}{2i\pi}\int_{\gamma}\vp(\omega)f(\omega)d\omega=\int_{S}\vp(x)\nu_{\beta, c'}(dx).$$
in particular, for $$\vp(\omega)=\frac{1}{z-\omega},$$
if $z$ is in the exterior of $\gamma$, then
$$\frac{1}{2i\pi}\int_{\gamma}\frac{1}{z-\omega}f(\omega)d\omega=G_{\beta, c'}(z).$$
We will use the theorem of residues to derive the expression of $G_{c'}$.

 The function $g(\omega)=\frac{1}{z-\omega}f(\omega)=\frac{2}{\beta\omega(z-\omega)}\sqrt{\omega-a}\sqrt{\omega-b}\sqrt{\omega+a}\sqrt{\omega+b},$ is meromorphic in $\overline{\Bbb C}\setminus S\cup\{0, z\} $. with simple pole at $\omega=0$, $\omega=z$ and a pole at infinity.

 Furthermore the residue at $\omega=0$ is $\displaystyle -\frac{2c}{\beta z}=-\frac{c'}{z}$, the residue at $\omega=z$ is $\displaystyle -f(z)$,
and $f$ admit a Laurent expansion for $|\omega|>\max(|a|,|b|,|z|)$
   $$f(\omega)=\frac{2\omega}{\beta}\sqrt{1-\frac{a}{\omega}}\sqrt{1-\frac{b}{\omega}}\sqrt{1+\frac{a}{\omega}}\sqrt{1+\frac{b}{\omega}}=
   \frac{2}{\beta}\omega-\frac{a^2+b^2}{\beta}\frac1\omega+\cdots$$
  then  the residue at $\omega=\infty$ is $\displaystyle-\frac{2z}{\beta}$.\\
  Which give $$G_{\beta, c'}(z)=-f(z)+\frac{2}{\beta}z-\frac{c'}{z}.$$

{\bf Proof of proposition 3.2.}

 Let denote by $U_{\beta, c'}$ the logarithmic potential of the measure $\nu_{\beta, c'}$: for all $x\in\Bbb R$,
 $$U_{\beta, c'}(x)=\int_\Bbb R\log\frac{1}{|x-t|}\nu_{\beta, c'}(dt).$$
 The function $U_{\beta, c'}$ is even and
 $$\frac{d}{dx}U_{\beta, c'}(x)=-\mbox {Re}G_{\beta, c'}(x).$$
 We will study the variation of the function
 $$\varphi(x)=U_{\beta, c'}(x)+\frac 12 Q_{c'}(x).$$
 The function $\varphi$ is even and
 $$\varphi'(x)=-\mbox {Re}G_{\beta, c'}(x)+\frac 12 \frac{d}{dx}Q_{c'}(x).$$
 (It is not defined on the point $x=0$). The last function vanished on $S$, therefore the function is constant on each connect components of $S$. Since the function $\varphi$ is even therefore the constant is the same on each components. Let denoted it by $C$.

 \begin{displaymath}
  \begin{array}{c|lllccccccr|}
    x    &\qquad -b&      &-a &            &0&           \quad &a             \qquad &b\qquad\\\hline
   \varphi'(x)     &\quad - &       0      &  &+            &  &-   & \quad 0 & & +   \\\hline
           &\quad\searrow&     &\quad\nearrow&   & &\quad\searrow&            &\quad\nearrow&                           \\
     \varphi(x)     &          &             C&          && &  & \quad C   &        \\
  \end{array}
\end{displaymath}
Therefore
 $$\begin{aligned} U_{\beta, c'}(x)+\frac 12 Q_{c'}(x)&\geq C\quad\mbox{in}\;\Bbb R,\\
&=C\quad\mbox{in}\; S.
\end{aligned}$$

By making use the proposition 2.2 the equilibrium measure $\nu^*$ coincide with $\nu_{\beta, c'}$.\qquad\qquad\qquad\qquad\qquad\qquad\qquad\qquad\qquad\qquad\qquad\qquad\qquad\qquad\qquad\qquad$\blacksquare$
\\
\subsection{Energy of equilibrium measure}
Consider the integral,
\small{\begin{equation}A_n=\int_{\Bbb R^n}e^{-n\sum\limits_{k=1}^{n}\lambda_k^2}\prod\limits_{k=1}^{n}|\lambda_k|^{2\mu_n}
\prod_{1\leq i<j\leq n}|\lambda_i-\lambda_j|^{\beta}d\lambda_1\cdots d\lambda_n=\int_{\Bbb R^n}e^{-K_n(\lambda)-\sum\limits_{i=1}^{n}Q_{ \alpha_n}(\lambda_i)}d\lambda_1\cdots\lambda_n.\end{equation}}
where
$$K_n(\lambda)=K_n(\lambda_1,\cdots,\lambda_n)=\frac{\beta}{2}\sum_{i\neq j}\log\frac{1}{|\lambda_i-\lambda_j|}+(n-1)\sum_{i=1}^nQ_{ \alpha_n}(\lambda_i),$$
and
 $$Q_{ \alpha_n}(x)= x^2+2\alpha_n\log\frac{1}{|x|},\qquad \alpha_n=\frac{\mu_n}{n}$$

 For $c\geq 0$ consider also the integral
 \begin{equation}B_n=\int_{\Bbb R^n}e^{-n\sum\limits_{k=1}^{n}\lambda_k^2}\prod\limits_{k=1}^{n}|\lambda_k|^{2nc}
\prod_{1\leq i<j\leq n}|\lambda_i-\lambda_j|^{\beta}d\lambda_1\cdots d\lambda_n=\int_{\Bbb R^n}e^{-n^2V_{n, c}(\lambda)}d\lambda_1\cdots\lambda_n.\end{equation}
where
 $$Q_{ c}(x)= x^2+2c\log\frac{1}{|x|},$$
Recall that the energy for a probability $\nu$ is defined by
$$E_{\beta, \delta}(\nu)=\frac{\beta}{2}\int_{\Bbb R^2}\log\frac{1}{|x-y|}\nu(dx)\nu(dy)+\int_{\Bbb R}Q_{\delta}(x)\nu(dx),$$
where  $$Q_{ \delta}(x)= x^2+2\delta\log\frac{1}{|x|},$$
We saw that $$\lim_{n\to\infty}-\frac{1}{n^2}\log B_n=E^*_{\beta, c}.$$
See for instance (Faraut \cite{F}). We will prove this result in proposition 4.6. for more general potential

Remark that $\displaystyle\lim_{x\to\pm\infty}K_n(x)=+\infty$  and $\displaystyle\lim_{x\to 0}K_n(x)=+\infty$, the same hold in the diagonal of $\Bbb R^n$.
Since the function $K_n$ is continuous except on the diagonal and $0$ where it has as limit $+\infty$. Hence it is bounded below and the minimum is realized at some point $\lambda^{(n)}=(\lambda^{(n)}_1,\cdots,\lambda^{(n)}_n)$, which means that
$$\inf_{\Bbb R^n}K_n(x)=K_n(\lambda^{(n)}).$$
Let denote by $$\tau_n=\frac{1}{n(n-1)}\inf_{x\in\Bbb R^n}K_n(x)\quad{\rm and}\quad\rho_n=\frac1n\sum_{i=1}^{n}\delta_{\lambda^{(n)}_i},$$
$\delta_{\lambda^{(n)}_i}$ is the Dirac mass at ${\lambda^{(n)}_i}$.\\

  From proposition 4.2, if we replace $c$ by $\alpha_n$ the equilibrium measure of the potential $\displaystyle \frac{2}{\beta}Q_{\alpha_n}$ is $\nu_{\beta, \alpha_n}$, where the density of the equilibrium measure $\nu_{\beta, \alpha_n}$ is given by
$$f_{\beta,\alpha_n}(t)=\left\{\begin{aligned}&\frac {2}{\pi\beta}\frac{1}{|t|}\sqrt{(t^2-a_n^2)(b_n^2-t^2)}\;\;\;if\;\;t\in S_n\\
&\quad 0\qquad\qquad\qquad\qquad\qquad\;\; if\;\;t\notin S_n\end{aligned}\right.,$$
 $S_n=[-b_n, a_n]\cup[a_n, b_n]$ and $a_n=\sqrt{\frac{\beta}{2}}\sqrt{1+\frac{\alpha_n}{\beta}-\sqrt {1+\frac{2\alpha_n}{\beta}}}$, $b_n=\sqrt{\frac{\beta}{2}}\sqrt{1+\frac{\alpha_n}{\beta}+\sqrt {1+\frac{2\alpha_n}{\beta}}}$.
\begin{lemma}
Let $(\mu_n)_n$ be a positive real sequence. Assume there is some constant $c$ such that $\displaystyle \lim_{n\to\infty}\frac{\mu_n}{n}=c.$ Then
\begin{description}
\item[(1)]  The probability measure $\nu_{\beta, \alpha_n}$ converge weakly to the probability $\nu_{\beta, c}$.
\item[(2)] $ E^*_{\beta,c}=\displaystyle\lim_{n\to\infty}E_{\beta, \alpha_n}.$

\end{description}
where $E^*_{\beta, c}$ is the energy of the equilibrium measure $\nu_{\beta, c}.$
\end{lemma}
\begin{proposition} Let $(\mu_n)_n$ be a positive real sequence. Assume there is some constant $c$ such that $\displaystyle \lim_{n\to\infty}\frac{\mu_n}{n}=c.$ Then
\begin{description}

\item[(1)] $\displaystyle \lim_{n\to\infty}\tau_n=E^*_{\beta, c}.$
\item[(2)] The measure $\rho_n$ converge weakly the the equilibrium measure $\nu_{\beta, c}$.
\item[(3)] $\displaystyle \displaystyle\lim_{n\to\infty}-\frac{1}{n^2}\log A_n=E^*_{\beta,c}.$
\end{description}
\end{proposition}
\begin{proposition} Let $(\mu_n)_n$ be a positive real sequence. Assume there is some constant $c$ such that $\displaystyle \lim_{n\to\infty}\frac{\mu_n}{n}=c.$ Then the energy $E^*_{\beta,c}$ is given by

$$\displaystyle E^*_{\beta,c}=\frac{3\beta}{8}+\frac{\beta}{4}\log(\frac{4}{\beta})+c(\frac{3}{2}+\log\frac4{\beta})+\frac{2c^2}\beta\log\frac{4c}{\beta}
-(\frac{2c^2}{\beta}+c+\frac\beta8)\log(1+\frac{4c}{\beta}).$$
\end{proposition}
For $c=0$, one recover's the energy of the $\beta$-Gaussian unitary ensemble
$$E^*_{\beta, 0}=\frac{3\beta}{8}+\displaystyle\frac{\beta}{4}\log(\frac{4}{\beta}).$$

{\bf Proof of lemma 3.5.}\\
${\bf Step(1):}$ The probability measures $\nu_{\beta, \alpha_n}$ and $\nu_{\beta, c}$ have density respectively $f_{\beta, \alpha_n}$ and $f_c$. It is easy to see that the density $f_{\beta, \alpha_n}$ converges Pointwise to the density $f_c$. Then by applying Fatou lemma we deduce the convergence in the weak topology.\\
${\bf Step (2):}$ We know by definition of the energy that
 \begin{equation}E^*_{\beta, c}=\inf_{\nu\in\mathfrak{M}^1(\Bbb R)}E_{\beta, c}(\nu)\leq E_{\beta, c}(\nu_{\beta, \alpha_n}),\end{equation}
and
$$E_{\beta, c}(\nu_{\beta, \alpha_n})=E_{\beta, \alpha_n}(\nu_{\beta, \alpha_n})+\int_{\Bbb R}\left(Q_{c}(x)-Q_{\alpha_n}\right)\nu_{\beta, \alpha_n}(dx),$$
 which can be writing as
 \begin{equation}E_{\beta, c}(\nu_{\beta, \alpha_n})=E^*_{\beta, \alpha_n}+\int_{\Bbb R}\left(Q_{c}(x)-Q_{\alpha_n}\right)\nu_{\beta, \alpha_n}(dx),\end{equation}
 where $\displaystyle E^*_{\beta, \alpha_n}=\inf_{\nu\in\mathfrak{M}^1(\Bbb R)}E_{\beta, \alpha_n}(\nu)=E_{\beta, \alpha_n}(\nu_{\beta, \alpha_n}).$

 Furthermore
\begin{equation}E^*_{\beta, \alpha_n}=\inf_{\nu\in\mathfrak{M}^1(\Bbb R)}E_{\beta, \alpha_n}(\nu)\leq E_{\beta, \alpha_n}(\nu_{\beta, c}),\end{equation}
and
$$E_{\beta, \alpha_n}(\nu_{\beta, c})=E_{\beta, c}(\nu_{\beta, c})+\int_{\Bbb R}\Big(Q_{\alpha_n}(x)-Q_{c}(x)\Big)\nu_{\beta, c}(dx),$$

\begin{equation}E_{\beta, \alpha_n}(\nu_{\beta, c})=E^*_{\beta, c}+\int_{\Bbb R}\left(Q_{\alpha_n}(x)-Q_{c}(x)\right)\nu_{\beta, c}(dx).\end{equation}
From equations $(3.3)$,  $(3.4)$,  $(3.5)$,  $(3.6)$ one gets
 \begin{equation}E^*_{\beta, c}+\int_{\Bbb R}\left(Q_{\alpha_n}(x)-Q_{c}(x)\right)\nu_{\beta, \alpha_n}(dx)\leq E^*_{\beta, \alpha_n}\leq E^*_{\beta, c}+\int_{\Bbb R}\left(Q_{\alpha_n}(x)-Q_{c}(x)\right)\nu_{\beta, c}(dx)\end{equation}
So it is enough to prove that the integrals go to $0$ when $n$ go to infinity.
Recall that the probability measures $\nu_{\beta, \alpha_n}$ and $\nu_{\beta, c}$ are supported respectively by $S_n$ and $S$.\\

Furthermore $$ |Q_{\alpha_n}(x)-Q_{c}(x)|=2|\alpha_n-c||\log|x||,$$
Since the sequence $b_n$ converge to $b$ hence there is some positive constant $C$ such that for all $n\in\Bbb N$,

$$\displaystyle \sup_{S_n\cup S}|\log|x||=\max(\log b_n,\log b)\leq C.$$
Take the limit in equation $(3.7)$ and use the facts that $\nu_{\beta, \alpha_n}$ and $\nu_{\beta, c}$ are probability measures and the sequence $\alpha_n$ converge to $c$ we deduce that
$$\lim_{n\to\infty}E^*_{\beta, \alpha_n}=E^*_{\beta, c}.$$
{\bf Proof of proposition 4.6.} \\
We will denote
$$k_\delta(s,t)=\log\frac{1}{|s-t|}+\frac12 Q_\delta(s)+\frac12 Q_\delta(t),$$
 for $\ell>0$,
$$k_\delta^{\ell}(s,t)=\inf(k_\delta(s,t), \ell).$$
and $$h_{\alpha_n}(t)=Q_{\alpha_n}(t)-\log(1+t^2).$$

{\bf Step(1):} In this step we will prove {\bf (1)} and {\bf (2)}.

Let $\gamma\in\mathfrak{M}^1(\Bbb R)$ be a probability measure then
$$\begin{aligned}\int_{\Bbb R^n}K_n(x)\gamma(dx_1)...\gamma(dx_n)&=\frac{\beta}{2}n(n-1)\int_{\Bbb R^2}\log\frac{1}{|s-t|}\gamma(ds)\gamma(dt)+n(n-1)\int_{\Bbb R}Q_{\alpha_n}(s)\gamma(ds)\\&=n(n-1)E_{\beta,\alpha_n}(\gamma).\end{aligned}$$
Then $$\tau_n\leq E_{\beta, \alpha_n}(\gamma),$$
for $\gamma=\nu_{\beta, \alpha_n}$  which is the equilibrium measure for the potential $Q_{\alpha_n}$, we obtains
$$\tau_n\leq E^*_{\beta, \alpha_n},$$
By using step 2 of lemma 3.5 we deduce
\begin{equation}\label{10}\limsup_{n}\tau_n\leq\lim E^*_{\beta, \alpha_n}=E^*_{\beta, c}.\end{equation}

Furthermore $$\begin{aligned}E^{\ell}_{\beta, \alpha_n}(\rho_n)&=\int_{\Bbb R^2}k^\ell_{\alpha_n}(s,t)\rho_n(ds)\rho_n(dt)\\
&=\frac1{n^2}\sum_{i,j=1}^{n}k_{\beta, \alpha_n}^\ell(\lambda_i^{(n)},\lambda_j^{(n)})\\
&\leq \frac1{n^2}\sum_{1\leq i\neq j\leq n}k_{\alpha_n}(\lambda_i^{(n)},\lambda_j^{(n)})+\frac{\ell}{n}\\
&=\frac1{n^2}K_n(\lambda^{(n)})+\frac{\ell}{n}\\
&=\frac{n(n-1)}{n^2}\tau_n+\frac{\ell}{n}\leq E^*_{\beta, \alpha_n}+\frac{\ell}{n} .
\end{aligned}$$
By the inequality
$$|s-t|\leq\sqrt{1+s^2}\sqrt{1+t^2},$$
it follows that \begin{equation}\frac12h_{\alpha_n}(s)+\frac12 h_{\alpha_n}(t)\leq k_{\alpha_n}(s,t),\end{equation}
and then
$$\int_{\Bbb R}h_{\alpha_n}(s)\rho_n(ds)\leq E^{\ell}_{\beta, \alpha_n}(\rho_n)\leq E^*_{\beta, \alpha_n}+\frac{\ell}{n}$$
Since the sequence  $E^*_{\beta, \alpha_n}+\frac{\ell}{n}$ is bounded uniformly on $n$ by some  positive constant $C_0$.
Furthermore $$h_{\alpha_n}(s)= Q_{\alpha_n}(s)-\log(1+s^2)=s^2+\alpha_n\log\frac{1}{|s|}-\log(1+s^2).$$

Since the positive sequence $\alpha_n$ converge to $c$, then there is two positive constants $a_1, a_2$ such that $a_1\leq \alpha_n\leq a_2$ and
$$h_{\alpha_n}(s)\geq s^2+a_1\log\frac{1}{|s|}-\log(1+s^2)=h_1(s)\quad\,{\rm if}\;0<|s|\leq 1$$
$$h_{\alpha_n}(s)\geq s^2+a_2\log\frac{1}{|s|}-\log(1+s^2)=h_2(s)\quad\,{\rm if}\;|s|\geq 1$$
Let $h(s)=\inf(h_1(s),h_2(s))$, then $\displaystyle\lim_{|s|\to\infty}h(s)=+\infty$ and
$$\int_{\Bbb R}h(s)\rho_n(ds)\leq C_0.$$
Hence by Prokhorov criterium this proves
that the sequence $(\rho_n)_n$ is relatively compact for the weak topology.
Therefore there is a converging subsequence:
$\rho_{n_k}$ to $\rho$ which means, for all bound continuous fonctions on $\Bbb R$
$$\lim_{n\to\infty}\int_{\Bbb R}f(x)\rho_{n_k}(dx)=\int_{\Bbb R}f(x)\rho(dx).$$
We will denote by $\rho_n$ the subsequence.

For $\ell>0$ consider as in the previous the kernel $\displaystyle k_{\alpha_n}^\ell(s,t)=\inf(k_{\alpha_n}(s,t),\ell)$ and $\displaystyle k_{c}^\ell(s,t)=\inf(k_c(s,t),\ell)$.\\
Let $\e>0$, there is $n_0$, such that for all $n\geq n_0$,
$$c-\e\leq\alpha_n\leq c+\e,$$
Let $n\geq n_0$, divided $\Sigma=\Bbb R^2\setminus\{(s,t)\mid s= t\;{\rm or }\;s=0\;{\rm or}\; t= 0\} $ to fourth region
$$R_1=\{(s,t)\in\Sigma\mid |s|\geq 1\,{\rm and}\; |t|\ge 1\},\quad R_2=\{(s,t)\in\Sigma\mid 0<|s|\leq 1\,{\rm and}\; 0<|t|\leq 1\},$$
and $$R_3=\{(s,t)\in\Sigma\mid 0<|s|\leq 1\,{\rm and}\; |t|\geq 1\},\quad R_4=\{(s,t)\in\Sigma\mid \,|s|\geq 1\,{\rm and}\; 0<|t|\leq 1\}.$$
If $(s,t)\in R_1$, then
$$k_{\alpha_n}(s, t)\geq k_{c+\e}(s,t).$$
If $(s,t)\in R_2$,
$$k_{\alpha_n}(s,t)\geq k_{c-\e}(s,t).$$
If $(s,t)\in R_3$,
$$k_{\alpha_n}(s,t)\geq \log\frac{1}{|s-t|}+\frac12Q_{c+\e}(t)+\frac 12Q_{c-\e}(s),$$
hence $$k_{\alpha_n}(s,t)\geq\frac12\left(k_{c+\e}(s,t)+k_{c-\e}(s,t)\right).$$
By symmetry of the kernel $k_{\alpha_n}$  the last inequality  is valid in $R_4$.\\
we obtain for $(s, t)\in\Sigma$,
$$k_{\alpha_n}(s, t)\geq  a\,k_{c+\e}(s,t)+b\,k_{c-\e}(s,t),$$
where $(a,b)=(1,0)$ in $R_1$, $(a,b)=(0,1)$ in $R_2$ and $(a,b)=(\frac12,\frac12)$ in $R_3\cup R_4$.
Hence if we take the infimum we obtain $$k^\ell_{\alpha_n}(s, t)\geq a\,k^\ell_{c+\e}(s,t)+b\,k^\ell_{c-\e}(s,t) .$$
Moreover for the energy one gets, for all $n\geq n_0$
$$a E^{\ell}_{\beta, c+\e}({\rho}_{n})+bE^{\ell}_{\beta, c-\e}({\rho}_{n})\leq E^{\ell}_{\beta, \alpha_{n}}({\rho}_{n}).$$
Which gives
$$a E^{\ell}_{\beta, c+\e}({\rho}_{n})+bE^{\ell}_{\beta, c-\e}({\rho}_{n})\leq\frac{n(n-1)}{n^2}\tau_n+\frac{\ell}{n}.$$
As $n$ goes to infinity we obtain
$$\liminf_n \Big(a E^{\ell}_{\beta, c+\e}({\rho}_{n})+bE^{\ell}_{\beta, c-\e}({\rho}_{n})\Big)\leq\liminf\tau_n,$$
hence by the weak convergence of the subsequence $\rho_n$ it follow
$$a E^{\ell}_{\beta, c+\e}({\rho})+bE^{\ell}_{\beta, c-\e}({\rho})\leq\liminf\tau_n,$$
applying the monotone convergence theorem, when $\ell$ goes to $0$, it follows that
$$a E_{\beta, c+\e}({\rho})+bE_{\beta, c-\e}({\rho})\leq\liminf\tau_n.$$
Since $\rho$ is a probability measure and using the values of $a, b$ we obtain $a E_{\beta, c+\e}({\rho})+bE_{\beta, c-\e}({\rho})=E_{\beta, c}({\rho}).$
hence $$E_{\beta, c}({\rho})\leq\liminf\tau_n.$$
Furthermore $$\inf_{\mu\in\mathfrak{M}^1(\Bbb R)}E_{\beta, c}(\mu)\leq E_{\beta, c}({\rho}).$$
We saw  from proposition 4.2. that the minimum is realized at the probability measure $\nu_{\beta, c}$ and the minimum is
$E^*_{\beta, c}.$
Hence $$E^*_{\beta, c}\leq E_{\beta, c}({\rho})\leq\liminf\tau_n.$$
It follows that $$E^*_{\beta, c}\leq E_{\beta, c}({\rho})\leq \liminf\tau_n\leq \limsup\tau_n\leq E^*_{\beta, c},$$
in the last inequalities we used equation (\ref{10}).
Therefore $$E({\rho})=E^*_{\beta, c}=E_{\beta, c}(\nu_{\beta, c}).$$
This implies that ${\rho}=\nu_{\beta, c}$. We have proved that $\nu_{\beta, c}$ is
the only possible limit for a subsequence of the sequence $(\rho_n)$. It follows
that the sequence $(\rho_n)$ itself converges:
for all bounded continuous  function $$\lim_{n\to\infty}\int_{\Bbb R}f(x){\rho}_n(dx)=\int_{\Bbb R}f(x)\nu_{\beta, c}(dx),$$
and $$\lim_{n\to\infty}\tau_n=E^*_{\beta, c}.$$

{\bf Step (2):} Now we will prove: $\displaystyle\lim_{n\to\infty}-\frac{1}{n^2}\log A_n=E^*_{\beta, c}$.

Recall that $$A_n=\int_{\Bbb R^n}e^{-K_n(\lambda)-\sum\limits_{i=1}^{n}Q_{\alpha_n}(\lambda_i)}\,d\lambda_1\cdots d\lambda_n,$$
it follows that
$$A_n\leq e^{-n(n-1)\tau_n}\left(\int_{\Bbb R}e^{-Q_{\alpha_n}(\lambda)}d\lambda)\right)^n=e^{-n(n-1)\tau_n}\left(\Gamma(\alpha_n+\frac12)\right)^n,$$
and
$$\frac{1}{n^2}\log A_n\leq -\frac{n-1}{n}\tau_n+\frac{1}{n}\log\Gamma(\alpha_n+\frac12).$$
Since the sequence $(\alpha_n)$ converge to $c$ then $\displaystyle\lim_{n\to\infty}\log\Gamma(\alpha_n+\frac12)=\Gamma(c+\frac12)$
which gives
\begin{equation}\label{112}\liminf_{n}-\frac{1}{n^2}\log A_n\geq\liminf_n\tau_n=E^*_{\beta, c}.\end{equation}
Furthermore if $\mu$ is a probability measure then
$$\int_{\Bbb R^n}K_n(x)\mu(dx_1)\cdots\mu(dx_n)=n(n-1)E_{\beta, \alpha_n}(\mu),$$
Let $\mu(dt)=\nu_{\beta, c}(dt)=f_{\beta, c}(t)dt$ supported by $S=[-b,-a]\cup[a,b]$, the function $f_{\beta, c}(t)>0$ except on subset of $S$ with measure zero.
Applying Jensen inequality to the exponential function then
$$\begin{aligned}A_n&=\int_{\Bbb R^n}\exp\left(-K_n(x)-\sum_{i=1}^{n}Q_{\alpha_n}(x_i)-\sum_{i=1}^{n}\log f_c(x_i)\right)\prod_{i=1}^nf_c(x_i)dx_1\cdots dx_n\\
&\geq\exp\left(\int_{\Bbb R^n}\left(-K_n(x)-\sum_{i=1}^{n}Q_{\alpha_n}(x_i)-\sum_{i=1}^{n}\log f_c(x_i)\right)\prod_{i=1}^nf_c(x_i)dx_1\cdots dx_n\right)\\
&\geq e^{-n(n-1)E_{\beta, \alpha_n}(\nu_{\beta, c})}\exp\left(-n\int_{\Bbb R}Q_{\alpha_n}(x)f_{\beta, c}(x)dx\right)\exp\left(-n\int_{\Bbb R}f_{\beta,c}(x)\log f_{\beta, c}(x)dx\right)
.\end{aligned}$$
From lemme 4.5 we have $$\lim_{n\to\infty}E_{\beta, \alpha_n}(\nu_{\beta, c})=E_{\beta, c}(\nu_{\beta, c})=E^*_{\beta, c},$$
and $$\int_{\Bbb R}Q_{\alpha_n}(x)f_c(x)dx=2\int_a^bQ_{\alpha_n}(x)f_c(x)dx\leq 2(b^2+2|\alpha_n|\log b),$$
furthermore
the last integral exist by the continuity of the function $x\log x$ near $0$ and the continuous function $f_c$ is with compactly support $S$. So the integral is bounded by some constant say $M$. Then
$$-\frac{1}{n^2}\log A_n\leq \frac{n-1}{n}E^*_{\beta, c}+ \frac1n\Big(2b^2+4|\alpha_n|\log b+M\Big).$$

It follows that
$$\limsup_{n}-\frac{1}{n^2}\log A_n\leq\limsup_{n}\Big( \frac{n-1}{n}E_{\beta, \alpha_n}(\nu_{\beta, c})+ \frac1n\big(b^2+|\alpha_n|\log a+M\big)\Big).$$
Since $\alpha_n$ converge. Hence
\begin{equation}\label{111}\limsup_{n}-\frac{1}{n^2}\log A_n\leq E^*_{\beta, c}.\end{equation}
Equations (\ref{112}) and (\ref{111}) gives that
$$\lim_{n}-\frac{1}{n^2}\log A_n= E^*_{\beta, c}.$$
Which complete the proof.\qquad\qquad\qquad\qquad\qquad\qquad\qquad\qquad\qquad\qquad\qquad\qquad\qquad\quad$\blacksquare$\\
If we choose $\mu_n=nc$ we obtains the same result for the sequence $B_n$,
$$\lim_{n}-\frac{1}{n^2}\log B_n= E^*_{\beta, c}.$$
{\bf Proof of proposition 4.7.}
For more convenient we will prove the proposition first when $\beta=2$ and then deduce from proposition 3.3. the result for all $\beta>0$ \\
{\bf First case $\beta=2$}. By performing the change of variable $x_k=\lambda_k\sqrt n$ in the expression of $A_n$ equation $(3.2)$, we obtain
$$A_n=n^{-n\mu_n-\frac{n^2}{2}}n!C_n=n^{-n\mu_n-\frac{n^2}{2}}n!\prod_{k=1}^{n-1}\gamma_{\mu_n}(k),$$
where $\gamma_{\mu_n}(k)$ is defined in equation $(3.1)$.\\
{\bf First step.} Let $n=2m$ be an even integer. Then
$$\begin{aligned}A_{2m}&=(2m)!(2m)^{-2m\mu_{2m}-\frac{(2m)^2}{2}}\prod_{k=1}^{2m-1}\gamma_{\mu_{2m}}(k)\\
&=(2m)!(2m)^{-2m\mu_{2m}-\frac{(2m)^2}{2}}\prod_{k=0}^{m-1}\gamma_{\mu_{2m}}(2k)\prod_{k=1}^{m-1}\gamma_{\mu_{2m}}(2k+1)\\
&=(2m)!(2m)^{-2m\mu_{2m}-\frac{(2m)^2}{2}}\prod_{k=0}^{m-1}(k!)^2(\Gamma(k+\mu_{2m}+\frac12))^2\prod_{k=1}^{m-1}(k+\mu_{2m}+\frac12),
\end{aligned}$$
in the last equality we use the fact that $\Gamma(x+1)=x\Gamma(x).$\\
Take the logarithm of $A_{2m}$
$$\begin{aligned}\log(A_{2m})&=\sum_{k=1}^{2m}\log(k)+2\sum_{k=1}^{m-1}(m-k)\log(k)+2\sum_{k=0}^{m-1}\log\Gamma(k+\mu_{2m}+\frac12)\\&+\sum_{k=0}^{m-1}
\log(k+\mu_{2m}+\frac12)
-(2m\mu_{2m}+\frac{(2m)^2}{2})\log(2m).
\end{aligned}$$
It is easy to see that for $m$ large enough
\begin{equation}\label{10}\sum_{k=1}^{2m}\log(k)=o(m^2).\end{equation}

Furthermore from the Stiriling asymptotic formula we have, for $0\leq k\leq m-1$
\begin{equation}\label{p}\log\Gamma(k+\mu_{2m}+\frac12)=(k+\mu_{2m})\log(k+\mu_{2m}+\frac12)
-(k+\mu_{2m}+\frac12)+o\left(\log(k+\mu_{2m}+\frac12)\right),\end{equation}
and by the fact that $\mu_n=cn+o(n)$, we deduce, that\\
$$ \log(k+\mu_{2m}+\frac12)=\log(k+\mu_{2m})+\log(1+\frac{1}{2(k+\mu_{2m})})=\log(k+\mu_{2m})+\frac{1}{2(k+\mu_{2m})}+o(\frac1m),$$
and $\displaystyle \log(k+\mu_{2m}+\frac12)=o(m),$  $\displaystyle\sum_{k=0}^{m-1}(k+\mu_{2m}+\frac12)=\sum_{k=0}^{m-1}(k+\mu_{2m})+o(m^2).$\\
By summing both side of (\ref{p}), one gets
\begin{equation}\label{101}\sum_{k=0}^{m-1}\log\Gamma(k+\mu_{2m}+\frac12)=\sum_{k=0}^{m-1}(k+\mu_{2m})\log(k+\mu_{2m})-\sum_{k=0}^{m-1}(k+\mu_{2m})+o(m^2),\end{equation}
and
\begin{equation}\label{102}\sum_{k=0}^{m-1}\log(k+\mu_{2m}+\frac12)=o(m^2).\end{equation}
Hence, from equation (\ref{10}), (\ref{101}) and  (\ref{102}), it follows
$$\begin{aligned}\frac{1}{(2m)^2}\log(A_{2m})=&-(\frac{\mu_{2m}}{2m}+\frac{1}{2})\log(2m)-\frac{2}{(2m)^2}\left(\frac{m(m-1)}{2}+m\mu_{2m}\right)\\
&+\frac{1}{2m}\sum_{k=1}^{m-1}(1-\frac{k}{m})\log(k)+\frac{1}{2m^2}\sum_{k=0}^{m-1}(k+\mu_{2m})\log(k+\mu_{2m})+o(1).
\end{aligned}$$
Thus
$$\begin{aligned}\frac{1}{(2m)^2}\log(A_{2m})=&-(\frac{\mu_{2m}}{2m}+\frac{1}{2})\log(2m)-(\frac12+\frac{\mu_{2m}}{2m})+\frac{1}{2m}\sum_{k=1}^{m-1}(1-\frac{k}{m})\log(m)\\
&+\frac{1}{2m^2}\sum_{k=0}^{m-1}(k+\mu_{2m})\log(m+\mu_{2m})+S_m^1+S_m^2+o(1),\end{aligned}$$
where
$$S_m^1=\frac{1}{2m}\sum_{k=1}^{m-1}(1-\frac{k}{m})\log(\frac{k}{m}),$$
and
$$S_m^2=\frac{1}{2m^2}\sum_{k=0}^{m-1}(k+\mu_{2m})\log(\frac{k+\mu_{2m}}{m+\mu_{2m}}).$$
Applying Riemann sums for both sums $S_m^1$ and $S_m^2$, we obtain
\begin{equation}\label{1}\lim_{m\to\infty}S_m^1=\lim_{m\to\infty}\frac{1}{2m}\sum_{k=1}^{m-1}(1-\frac{k}{m})\log(\frac{k}{m})=\frac12\int_0^1(1-x)\log xdx=-\frac38,\end{equation}
\begin{equation}\label{2}\begin{aligned}\lim_{m\to\infty}S_m^2&=\lim_{m\to\infty}\frac12\left(1+\frac{\mu_{2m}}{m}\right)^2\frac{1}{m+\mu_{2m}}\sum_{k=0}^{m-1}(\frac{k+\mu_{2m}}{m+\mu_{2m}})\log(\frac{k+\mu_{2m}}{m+\mu_{2m}})\\
&=\frac12(1+2c)^2\int_{\frac{2c}{1+2c}}^1x\log xdx=-\frac18(1+2c)^2+\frac12c^2+c^2\log(1+\frac1{2c}).\end{aligned}\end{equation}
Now we will compute the limits of the others terms
$$I_m=-(\frac{\mu_{2m}}{2m}+\frac{1}{2})\log(2m)+\frac{1}{2m}\sum_{k=1}^{m-1}(1-\frac{k}{m})\log(m)
+\frac{1}{2m^2}\sum_{k=0}^{m-1}(k+\mu_{2m})\log(m+\mu_{2m})-(\frac14+\frac{\mu_{2m}}{2m})$$
By simple computation it  yields
$$I_m=-(\frac{\mu_{2m}}{2m}+\frac{1}{2})\log(2m)+\frac{m-1}{4m}\log(m)+\frac{1}{2m^2}\left(\frac{m(m-1)}{2}
+m\mu_{2m}\right)\log(m+\mu_{2m})-(\frac14+\frac{\mu_{2m}}{2m}).$$
$$\begin{aligned}I_m=-(\frac{\mu_{2m}}{2m}+\frac{1}{2})\log2+\left(\frac{m-1}{2m}-\frac12\right)\log(m)
+\left(\frac{m-1}{4m}+\frac{\mu_{2m}}{2m}\right)\log(1+\frac{\mu_{2m}}{m})-\frac14-\frac{\mu_{2m}}{2m}.\end{aligned}$$
Hence
\begin{equation}\label{3}\lim_{m\to\infty}I_m=-(c+\frac{1}{2})\log2+(\frac14+c)\log(1+2c)-\frac14-c.\end{equation}

From equations (\ref{1}), (\ref{2}) and (\ref{3}) it follows
$$\lim_{m\to\infty}-\frac{1}{(2m)^2}\log A_{2m}=\frac34+\frac12\log2+(\frac{3}{2}+\log2)c+c^2\log(2c)-(c^2+c+\frac14)\log(1+2c).$$
\\
{\bf Second step.} when $n=2m+1$, we prove by the same method that $$\displaystyle \lim_{m\to\infty}-\frac{1}{(2m+1)^2}\log A_{2m+1}=\frac34+\frac12\log2+(\frac{3}{2}+\log2)c+c^2\log(2c)-(c^2+c+\frac14)\log(1+2c).$$
Furthermore it is easy to see that the integral $B_n$ is a particular case of $A_n$ when we take $\mu_n=nc$. Then we have
$$\lim_{n\to\infty}-\frac{1}{n^2}\log B_n=\lim_{n\to\infty}-\frac{1}{n^2}\log A_n=E^*_{2, c}.$$

{\bf Second case $\beta>0$.}
 Define for $\nu\in\mathfrak{M}^1(\Bbb R)$ the energy
$$E_{\beta, \alpha_n}(\nu)=\frac{\beta}{2}\left(\int_{\Bbb R^2}\log\frac{1}{|s-t|}\nu(ds)\nu(dt)+\int_\Bbb RQ_{\beta, \alpha_n}(t)\nu(dt)\right),$$
where $$\displaystyle Q_{\beta,\alpha_n}(t)=\left(\sqrt{\frac{2}{\beta}}t\right)^2+\frac{4\alpha_n}{\beta}\log\frac{1}{|t|}.$$\\
Since
$$\displaystyle Q_{\beta, \alpha_n}(t)=Q_{2,\alpha_n}\circ h(t)+\frac{4\alpha_n}{\beta }\log\sqrt{\frac 2\beta},$$
where $h(t)=\sqrt{\frac{2}{\beta}}t$. Then
by proposition 3.3, we obtains
$$E_{\beta, \alpha_n }=\frac{\beta}{2}E_{2,\frac{2\alpha_n}{\beta}}+\frac{\beta}{2}\log\sqrt{\frac{2}{\beta}}+2\alpha_n\log\sqrt{\frac 2\beta}.$$
We saw from lemma 3.5  that
$$\lim_{n\to\infty}E_{2,\frac{2\alpha_n}{\beta}}=E^*_{2,\frac{2c}{\beta}}\; {\rm and}\; \lim_{n\to\infty}E_{\beta,\alpha_n}=E^*_{\beta,c}.$$
From the first case $\beta=2$ and simple computation we deduce the desired result.$\blacksquare$
\section{Proof of theorem 4.1}
 Recall the statistical
 distribution $\nu_n$ is defined by: for all bounded continuous  function $f$ on $\Bbb R$,
 $$\int_{\Bbb R}f(t)\nu_n(dt)=\mathbb{E}_{n,\mu_n}\left(\frac1n\sum_{i=1}^nf(\lambda_i)\right),$$
 where $\mathbb{E}_{n,\mu_n}$ is the expectation with respect the probability on $\Bbb R^n$

$$\mathbb{P}_{n,\mu_n}(d\lambda)=\frac{1}{Z_n}e^{-n\sum\limits_{i=1}^{n}\lambda_i^2}\prod\limits_{i=1}^{n}|\lambda_i|^{2\mu_n}
\prod_{1\leq i<j\leq n}|\lambda_i-\lambda_j|^{\beta}d\lambda_1\cdots d\lambda_n.$$
Let Define on $\Bbb R^n$ the function :
$$K_n(x)=\frac{\beta}{2}\sum_{i\neq j}\log\frac{1}{|x_i-x_j|}+(n-1)\sum_{i=1}^{n}Q_{\alpha_n}(x_i),$$
where $\displaystyle Q_{\alpha_n}=x^2+ 2\alpha_n\log\frac{1}{|x|}$
and $\displaystyle \alpha_n=\frac{\mu_n}{n}.$\\
The probability $\mathbb{P}_{n,\mu_n}$ concentrates in a neighborhood of the points where
the function $\displaystyle K_n(x)$
attains its infimum:
\begin{proposition} Let $\e>0$ and $\displaystyle A_{n,\e}=\left\{x\in\Bbb R^n\mid K_n(x)\leq (E^*_{\beta, c}+\e)n^2\right\}$. Then\\
$A_{n,\e}$ is compact and $$\lim_{n\to\infty}\mathbb{P}_{n,\mu_n}(A_{n,\e})=1.$$
\end{proposition}
This proposition can be found in \cite{F}, lemma IV.5.2. We give the proof.\\
{\bf Proof.}
Recall that $h_{\alpha_n}(x)=Q_{\alpha_n}(x)-\log(1+x^2)$. Since $h_{\alpha_n}$ is lower semicontinuous and
$$K_n(x)\geq(n-1)\sum_{i=1}^n h_{\alpha_n}(x_i),\qquad \lim_{x_i\to\pm\infty}h_{\alpha_n}(x_i)=+\infty,$$
 then $A_{n, \e}$ is closed and bounded hence it is compact.

Let $\e>0$, from the definition of $A_{n,\e}$ we have on $\Bbb R^n\setminus A_{n,\e}$
$$K_n(x)>(E^*_{\beta, c}+\e)n^2,$$ then
$$\mathbb{P}_{n,\mu_n}(\Bbb R^n\setminus A_{n,\e})\leq \frac{1}{Z_n}e^{-(E^*_{\beta, c}+\e)n^2}\left(\int_{\Bbb R}e^{-Q_{\alpha_n}(x)dx}\right)^n.$$
Furthermore
$$\int_{\Bbb R}e^{-Q_{\alpha_n}(x)}dx=\int_{\Bbb R}|x|^{2\alpha_n}e^{-x^2}dx=\Gamma(\alpha_n+\frac12).$$
By continuity of the gamma function we have $\displaystyle\lim_{n\to\infty}\Gamma(\alpha_n+\frac12)=\Gamma(c+\frac12).$
Since from proposition 4.6  $\displaystyle\lim_{n\to\infty}-\frac{1}{n^2}\log Z_n=E^*_{\beta,c}$. Then there is $n_0$ such that for all $n\geq n_0$,
$$\frac{1}{Z_n}\leq e^{\big(E^*_{\beta, c}+\frac{\e}{2}\big)n^2}.$$
Using all those arguments we obtain for $n$ large enough
$$\mathbb{P}_{n,\mu_n}(\Bbb R^n\setminus A_{n,\e})\leq\left(\Gamma(c+\frac12)+\e\right)^ne^{-\frac{\e}{2}n^2}.$$
Which complete the proof.\qquad\qquad\qquad\qquad\qquad\qquad\qquad\qquad\qquad\qquad\qquad\qquad\qquad\qquad\quad$\blacksquare$\\\\
{\bf Proof of theorem 4.1.}
We keep those notations:

 $$k_{\alpha_n}(s,t)=\log\frac{1}{|s-t|}+\frac12 Q_{\alpha_n}(s)+\frac12 Q_{\alpha_n}(t),$$
 for $\ell>0$,
$$k_{\alpha_n}^{\ell}(s,t)=\inf(k_\delta(s,t), \ell).$$
 $$h_{\alpha_n}(t)=Q_{\alpha_n}(t)-\log(1+t^2),$$
and $h(t)=\inf(h_{a_1}(t),h_{a_2}(t))$, $h_{a_1}$ and $h_{a_2}$ are the functions used on the proof of proposition 3.6 where $a_1, a_2\geq 0$\\

For a bounded continuous function $f$ on $\Bbb R$, defined on $\Bbb R^n$ the continuous function
$$F_n(x)=\frac1n\sum_{i=1}^nf(x_i),$$
Let $\e>$, the set $A_{n, \e}$ is compact, hence $F_n$ attaint it supremum at same point in $A_{n, \e}$ say
$$x^{(n)}_{\e}=(x^{(n)}_{1, \e},\cdots,x^{(n)}_{n, \e}).$$
We obtain
$$\int_{\Bbb R}f(t)\nu_n(dt)\leq F_n(x^{(n)}_{\e})+||f||_{\infty}(1-\mathbb{P}_{n, \mu_n}(A_{n, \e})).$$
To the point $x^{(n)}_{\e}$ we associate the probability measure on $\Bbb R$
$$\sigma_{n, \e}=\frac1n\sum_{i=1}^n\delta_{ x^{(n)}_{i,\e}}.$$
The previous inequality can be written
$$\int_{\Bbb R}f(t)\nu_n(dt)\leq\int_{\Bbb R}f(t)\sigma_{n, \e}(dt)+||f||_{\infty}(1-\mathbb{P}_{n, \mu_n}(A_{n, \e})),$$
The truncated energy $E^{\ell}$ of the measure $\sigma_{n, \e}$ satisfies:
$$E^{\ell}(\sigma_{n, \e})\leq\frac{\ell}{n}+(E^*_{\beta, c}+\e).$$
From the inequality
$$(n-1)\sum_{i=1}^nh(x_i)\leq K_n(x),$$
we obtain
$$\int_{\Bbb R}h(t)\sigma_{n, \e}(dt)\leq \frac{n}{n-1}(E^*_{\beta, c}+\e).$$
This implies that the sequence $\sigma_{n, \e}$ is relatively compact for the weak topology. There is a sequence $n_j$ going to $\infty$ such that the subsequence $\sigma_{n_j, \e}$ converges in the weak topology: $$\lim_{n\to\infty}\sigma_{n_j, \e} =\sigma_\e.$$
We may also assume  in the weak topology that
$$\lim_{j\to\infty}\nu_{n_j}=\limsup_n\nu_n.$$
the limit measure satisfies
$$E^{\ell}(\sigma_{\e})\leq E^*_{\beta, c}+\e.$$

The rest of the proof is analogous to the proof of theorem IV.5.1 \cite{F}.

\begin{tabular}{ccc}

  College of Applied Sciences &&\qquad\qquad\qquad\qquad\qquad Faculte des Sciences de Tunis \\
   Umm Al-Qura University&&\qquad\qquad\qquad\qquad\qquad  Universit\'e de Tunis\\
  P.O Box  (715), Makkah &&\qquad\qquad\qquad\qquad\qquad  Campus Universitaire 2092, \\
  Saudi Arabia && \qquad\qquad\qquad\qquad\qquad   El Manar Tunis \\
  mabouali@uqu.edu.sa && \qquad\qquad\qquad\qquad\qquad  bouali25@laposte.net \\

\end{tabular}


\begin{thebibliography}{9}
\bibitem{I0} M. Bouali, \emph{Generalized Gaussian Random Unitary Matrices Ensemble}, (2014). To appear in JP Journal of Geometry and Topology.
\bibitem{II}
I. Dumitriu, \emph{Eigenvalue Statistics for Beta-Ensembles}, Ph.D thesis, (2003) Department of Mathematics, MIT.
\bibitem{I}
 I. Dumitriu, \& A. Edelman, \emph{ Matrix models for beta ensembles}, J. Math. Phys. 43, (2002), 5830-5847.
\bibitem{F}
J. Faraut, \emph{ Logarithmic potential theory, orthogonal polynomials,
and random matrices}, CIMPA Scool, Hammamet (2011).

\bibitem{fe} W. Feller, \emph{An introduction to probability theory and its applications,} Vol. II, Wisely (1971).

\bibitem{fo}P. J. Forrester, P. J, \emph{Log-gaes and random matrices}. (LMS-34). London Mathematical Sociaty monographs. Princeton University Press (2010).
\bibitem{H}
     U. Haagerup \& S. Thorbjørnsen,
      \emph{Random matrices with
complex Gaussian entries}, Exp. Mat. 21, (2003), 293-337.
\bibitem{K}
K. Johansson, \emph{On fluctuation of eigenvalues of random Hermitian matrices}, Duke. Math. J 91, (1998), 151-204.


\bibitem{M}
 M. L. Mehta, \emph{Random matrices}, Second Edition, Academic Press, (1991).


\bibitem{R}
M. Rosenblum, \emph{Generalized Hermite Polynomials and the Bose-Like Oscillator Calculus},
 Nonselfadjoint Operators and Related Topics
Operator Theory: Advances and Applications Volume 73, pp 369-396, (1994).

\bibitem{sa} E. B. Saff, E.B. \& V. Totik,  Logarithmic potentials with external
fields, Springer, (1997).





\end{thebibliography}
\end{document}